\documentclass[a4paper, 12pt]{amsart}
\usepackage{amssymb,amsmath,amsthm,amsfonts,latexsym}
\usepackage[utf8x]{inputenc}
\usepackage{mathrsfs}

\usepackage{setspace}

\usepackage[hyphens]{url}

\usepackage{enumitem}

\newtheorem{definition}{Definition}[section]
\newtheorem{lemma}{Lemma}[section]
\newtheorem{claim}{Claim}[section]
\newtheorem{proposition}{Proposition}[section]
\newtheorem{remark}{Remark}[section]

\newtheorem{corollary}{Corollary}[section]

\newcommand{\union}[1]{\bigcup\limits_{#1}}

\newtheorem{theorem}{Theorem}[section]

\title{The P-Ideal Dichotomy and the Failure of Two-Cardinal Square Principles}
\author{V\'ictor Torres-P\'erez}
\author{Liuzhen Wu}

\begin{document}


\begin{abstract}
We prove that the $\mathrm{P}$-ideal dichotomy implies the  negation of the square principle  $\square(\lambda, \omega)$ for every regular cardinal $\lambda\geq\omega_2$.
\end{abstract}

\maketitle

\section{Introduction}

The P-ideal dichotomy (PID) was introduced by  Todor{\v{c}}evi{\'c} in \cite{todorvcevic2000dichotomy}.  PID is a consequence of PFA (\cite{todorvcevic2000dichotomy}), and it retains many of its important consequences: it implies the negation of $\square(\lambda)$ for every regular cardinal $\lambda\geq\omega_2$ (\cite{todorvcevic2000dichotomy}), the Singular Cardinal Hypothesis \cite{viale2008family}, that the bounding number is at most $\aleph_2$ (see for example Corollary 21.8 in \cite{todorcevic2011combinatorial}), etc. 

One characteristic of PID is that it is compatible with the Continuum Hypothesis (\cite{todorvcevic2000dichotomy}).   This gives this principle considerable flexibility in contrast to traditional forcing axioms.  Since $\mathrm{CH}$ implies the negation of $\mathrm{MA}_{\aleph_1}$,  it is of interest in  the context of our previous results, where we worked with principles strong enough to have similar consequences to traditional forcing axioms like MM and PFA, but they can coexist with the negation of $\mathrm{MA}_{\aleph_1}$.  

Previously, we have worked, together with other co-authors,   with Rado's Conjecture (RC) or one of its main consequences, a strong version of Chang's Conjecture ($\mathrm{CC}^*$). For example, similar to MM, it was possible to establish a general panorama of the relationship of  RC with Schimmerling's square principle $\square_{\kappa, \lambda}$ (see \cite{sakai}).  We remark that $\mathrm{RC}$ implies the negation of $\mathrm{MA}_{\omega_1}$ (see  p. 8  in \cite{stevoradoconsistent}).  

Regarding sequences of the form $\square(\lambda, \mu)$, we  proved that $\mathrm{RC}$ implies $\lnot\square(\lambda,\omega)$ for every regular cardinal $\lambda\geq\omega_2$ (\cite{torreswu2}). 

Returning to PID,  
Raghavan proved among other results that $\mathrm{PID}$ implies  $\lnot\square_{\kappa, \omega}$ for all uncountable cardinals $\kappa$ (\cite{raghavan}). Since $\square_{\kappa, \mu}$ implies $\square(\kappa^+, \mu)$,  a natural question is whether $\mathrm{PID}$ proves the stronger result $\lnot\square(\lambda, \omega)$ for all regular $\lambda\geq\omega_2$ (Question 19.b in \cite{raghavan}) . We prove in this paper that in fact this is the case.  

\section{Definitions, notations and basic properties}

\subsection{Square sequences}
We recall a two-cardinal version of the square principle of the form $\square(\lambda, \mu)$.

\begin{definition}\label{defsquare}
Given a regular cardinal $\lambda$, a cardinal $\mu\leq\lambda$, 
 $\langle\mathscr C_\alpha: \alpha\in \mathrm{Lim}(\lambda)\rangle$ is a \emph{$\left(\lambda, \mu\right)$ $C$-sequence}  if 
\begin{enumerate}
\item $1\leq | \mathscr C_\alpha |\leq\mu$,
\item  for every $C\in \mathscr C_\alpha$, $C$ is a closed and unbounded subset of $\alpha$,
\item for every $C\in \mathscr C_\beta$, if $\alpha\in \mathrm{Lim}(C)$, then  $C\cap\alpha\in \mathscr C_\alpha$.
\end{enumerate}

For the sake of completeness, we also  define $\mathscr C_{\alpha+1}= \{\{\alpha\}\}$ for every $\alpha\in \lambda$. 

Given a set $C\subseteq \lambda$, we say $C$ \emph{trivializes}  a  $(\lambda, \mu)$ $C$-sequence
$\langle \mathscr C_\alpha: \alpha\in \mathrm{Lim}(\lambda)\rangle$   if $C\cap \alpha\in \mathscr C_\alpha$ for every $\alpha\in  \mathrm{Lim}(C).$

We say that the principle $\square(\lambda, \mu)$ holds if there is a $\left(\lambda, \mu\right)$ $C$-sequence $\langle \mathscr C_\alpha: \alpha\in \mathrm{Lim}(\lambda)\rangle$  (a $\square(\lambda, \mu)$-\emph{sequence}) such that there is no club trivializing this sequence. 
\end{definition}

\subsection{P-ideal dichotomy}

\begin{definition}
For two sets $x, y$, we write $x\subseteq^* y $ if $x\backslash y$ is finite. 
\end{definition}

\begin{definition}
 Let $X$ be an uncountable set. An ideal $\mathcal I\subseteq \left[ X\right]^{\leq\omega}$ is called a \emph{$P$-ideal} if for every collection $\left\{x_n: n\in \omega\right\}\subseteq \mathcal I$, there is $x\in \mathcal I$
 such that $x_n\subseteq^* x$ for all $n\in \omega$. 
\end{definition}

All ideals considered here contain every finite subset of their underlying
set,  
i.e. $\left[X\right]^{<\omega}\subseteq \mathcal I$.
\begin{definition}
 The P-ideal dichotomy $(\mathrm{PID})$ is the following statement: For any $P$-ideal $\mathcal I$ on an uncountable set $X$, either
 \begin{enumerate}
  \item there is an uncountable set $Y\subseteq X$ such that $\left[Y\right]^{\leq \omega}\subseteq \mathcal I$, or
  \item there exists  a partition $\{X_n:n\in \omega\}$ of $X$ such that $\left[X_n\right]^\omega\cap \mathcal I=\emptyset$ for all $n\in \omega.$
 \end{enumerate}

\end{definition}

In this paper, we prove the following: 

\begin{theorem}\label{pidnonsquare}
$\mathrm{PID}$ implies $\lnot\square(\lambda, \omega)$ for all regular cardinals $\lambda\geq \omega_2$. 
\end{theorem}

\section{PID and square sequences}

Fix a $(\lambda, \mu)$ $C$-sequence $\langle\mathscr C_\alpha: \alpha<\lambda\rangle$  and for every $\alpha\in \lambda$, let $\mathscr C_\alpha= \{C^\xi_\alpha: \xi\in \mu\}$.   Denote by $\mathrm{FIN}$ the set $\left[\mu\right]^{<\omega}\backslash\{\emptyset\}$. 

\begin{definition}
 For $F\in \mathrm{FIN}$ and $\alpha\leq\beta< \lambda$, define $T_F(\alpha, \beta)$ as follows: $T_F(\alpha, \alpha)=\{\langle \alpha\rangle\}$.
 For $\alpha<\beta$, let 
 $$T_F(\alpha, \beta)= 
	      \left\{
		      \langle \beta\rangle^\frown\sigma: \exists i\in F
			  \left(
			  \sigma \in 
				    T_F
				    \left(
					  \alpha, 
					  \min
					      \left(
						    C^i_\beta\backslash \alpha
					      \right)
				    \right)
			  \right)
	      \right\}.
	      $$
	      
\end{definition}

\begin{definition}
 For $F\in \mathrm{FIN}$ and $\alpha\leq\beta<\lambda$, let
 $$
 S_F(\alpha, \beta)= \min\{ |\sigma|: \sigma\in T_F(\alpha,\beta)\}.
 $$
\end{definition}

\begin{definition}
For $X\in [\lambda]^{\leq \omega}$, $\beta\in \lambda$, $F\in \mathrm{FIN}$ and $k\in \omega$, let 
$$
 H(X,  \beta, F,  k):=_{def}\{\alpha\in X\cap (\beta+1): S_F(\alpha, \beta)= k\}. 
 $$
\end{definition}

\begin{definition}

 Define $\mathcal I$ to be the collection of sets $X\in [\lambda]^{\leq \omega}$ such that  the set 
 $H(X, \beta, F,  k)$
 is finite for all 
 $\langle \beta,  F,  k\rangle\in  \mathrm{Lim}(X)\times \mathrm{FIN}\times   \omega $. 

\end{definition}

The proof of the following lemma does not depend on the order
types of the clubs, the cardinalities of the families of clubs,
or the cardinality of $\mathrm{FIN}$.

\begin{lemma}\label{lemma1}
Fix $X\in \mathcal I$ and $\gamma<\lambda$.  If $\sup(X)\leq \gamma$, then for all $F\in \mathrm{FIN}$ and for all $k\in\omega$, the set $\{\alpha\in X: S_F(\alpha, \gamma)= k\}$ is finite. 
\end{lemma}

\begin{proof}
See the  proof of Lemma 10 in \cite{raghavan}. 
\end{proof}

\begin{remark}\label{remark1}
$\mathcal I $ is closed under subsets. \footnote{Suppose   $Y\subseteq X$ with $X\in \mathcal I $.  Then $\mathrm{Lim}(Y)\subseteq \mathrm{Lim}(X)$, and so $H(Y, \beta, F, k)\subseteq H(X, \beta, F, k)$ for all $\langle \beta, F, k\rangle \in  \mathrm{Lim}(Y)\times \mathrm{FIN}\times\omega $.  
Since  $H(X, \beta, F, k)$ is finite,  we have $H(Y, \beta, F, k)$ is also finite, and therefore  $Y\in\mathcal I  $.  }
\end{remark}

Corollary \ref{corollary}, Lemma \ref{lemmaideal} and Lemma \ref{dilip} are also in \cite{raghavan}. We include their proofs for the benefit of the reader.

The following corollary is crucial, as it is used in many proofs related to these results. It demonstrates that the set $H(X, \gamma, F, k)$ is finite, even if $\gamma$ is not a limit point of $X,$ as long as $X$ is in the ideal $\mathcal{I}$. We also note that neither the order type of the clubs nor the size of $\mathrm{FIN}$ is taken into account.

\begin{corollary} \label{corollary}
The set   
$$H(X, \gamma, F, k)$$ 
is finite for every $\langle X, \gamma, F, k\rangle\in \mathcal I\times \lambda\times \mathrm{FIN}\times\omega$.  
\end{corollary}

\begin{proof}
Suppose for a contradiction that $H(X, \gamma, F, k)$ is infinite. Since $H(X, \gamma, F, k)\subseteq X$, then  by Remark  \ref{remark1}, $H(X, \gamma, F, k)\in \mathcal I$, contradicting Lemma \ref{lemma1}. 
 
\end{proof}

\begin{lemma}[Raghavan]\label{lemmaideal}
$\mathcal I$ is an ideal. 
\end{lemma}

\begin{proof}
Remark \ref{remark1} shows that  $\mathcal I$ is closed under subsets.  
We show now that it is closed under the union operator.  Take $X, Y\in \mathcal I$. Let $\beta\in \mathrm{Lim}(X\cup Y)$,  $F\in \mathrm{FIN}$  and $k\in \omega$. By Corollary \ref{corollary}, $H(X, \beta, F, k)$ and   $H(Y, \beta, F, k)$ are finite. 
Since $H(X\cup Y, \beta, F, k)\subseteq H(X, \beta, F, k)\cup H(Y, \beta, F, k)$, we have $X\cup Y\in \mathcal I$. 
\end{proof}

We recall the concept of the \emph{bounding number}. For two functions $f, g: \omega \to \omega$, we write $f \leq^* g$ if the set $\{n \in \omega : f(n) > g(n)\}$ is finite. The bounding number $\mathfrak{b}$ is defined as
$$
\mathfrak{b} = \min\left\{|\mathcal{F}| : \mathcal{F} \subseteq \omega^\omega \text{ and } \forall f \in \omega^\omega, \exists g \in \mathcal{F} \text{ such that } g \not\leq^* f\right\}.
$$

The proofs presented by Raghavan to show that $\mathcal I $ is a P-ideal in different instances (Lemma 12 and Theorem 20 in \cite{raghavan}) can be merged into a single, general proof. 

\begin{lemma}[Raghavan]\label{dilip}
Assume  $\mu<\mathfrak b$.  Then $\mathcal I$ is a $P$-ideal. 
\end{lemma}

\begin{proof}

Fix $\{X_n : n \in \omega\} \subseteq \mathcal{I}$. We have to find  $Z \in \mathcal{I}$ such that $X_n \subseteq^* Z$ for all $n \in \omega$. Without loss of generality, we may assume that $\{X_n: n\in \omega\}$ is  a collection of infinite sets.

Let $$Y= \union {n\in \omega}X_n,$$  
$$\mathrm{Lim}(Y)= \{\beta_m:m\in \omega\},$$
and 
for $n\in\omega$, 
$$X_n= \{\alpha^l_n:l\in\omega\}.$$

 By Corollary \ref{corollary}, for every $\langle n, m, F, k\rangle\in \omega\times\omega\times \mathrm{FIN}\times \omega$, the set   $H(X_n, \beta_m, F, k)$ is finite, 
and so we can build functions $f_{m, F, k}:\omega\to\omega$ such that   $H(X_n, \beta_m, F, k)\subseteq \{\alpha^l_n: l<f_{m, F, k}(n)\}$, for every $\langle n, m, F, k\rangle\in \omega\times\omega\times\mathrm{FIN}\times\omega$.
\footnote{
For example, define 
$f_{m,F,k}:\omega\to\omega$ 
by
$f_{m,F,k}(n)=
\max(\{0\}\cup
\{l+1:\alpha^l_n\in H(X_n,\beta_m,F,k)\})$.
} 
Since $\mu<\mathfrak b$,  the set $\{ f_{m, F, k}:  \langle m, F, k\rangle\in \omega\times\mathrm{FIN}\times\omega\}$ is a collection of $<\mathfrak b$ functions from $\omega$ to $\omega$, and  therefore we can take $f: \omega\to\omega$ such that $f_{m, F, k}\leq^*f$ for every $\langle m, F, k\rangle \in \omega\times \mathrm{FIN}\times \omega$. 

Let 
$$H(n)=\{\alpha^l_n: l<f(n)\} $$ 
and make 

$$Z=\union{n\in\omega}\left(X_n\backslash H(n)\right).$$

Observe that $X_n\setminus Z\subseteq H(n)$ for every
$n\in\omega$. Since $H(n)$ is finite, $X_n\subseteq^*Z$.

We  now show  that $Z$   is a member of $\mathcal I$. Let $(F, \beta, k)\in \mathrm{FIN}\times \mathrm{Lim}(Z)\times \omega. $
We must verify  that $H(Z,  \beta, F, k)$ is finite.

Observe that since $Z\subseteq Y$, then $\mathrm{Lim}(Z)\subseteq \mathrm{Lim}(Y)$, and so there is $m\in \omega$ such that $\beta= \beta_m$. 

\begin{claim}
    There is $N\in \omega$ such that  $H(X_n, \beta_m, F, k) \subseteq H(n) $ for every $n\geq N$. 
\end{claim}
\begin{proof}
Since $f_{m, F, k}\leq^* f$, take $N\in \omega$ such that $f_{m, F, k}(n)\leq f(n)$ for every $n\geq N$. Then if $\alpha^l_n\in H(X_n, \beta_m, F, k) $, we have $l<f_{m, F, k}(n)\leq f(n)$,  by the  definition of $f_{m, F, k}$ and  the construction  
 of $f$. Therefore, $\alpha^l_n\in H(n)$ by  the definition of $H(n)$.   
\end{proof}

By the preceding Claim,  take $N\in \omega$ such that $H(X_n, \beta_m, F, k)\subseteq  H(n) $ for every $ n\geq N$.  Then 
$ H(Z, \beta,F,  k)
    \subseteq 
    \union{n<N}
    H(X_n,  \beta_m,F,  k),$  
which is a finite union of finite sets. Therefore 
$ H(Z, \beta, F,  k)$ is finite and so $Z\in \mathcal I$. This finishes the proof of the Lemma. 

\end{proof}

The following lemma rules out the first alternative of PID.
Its proof is exactly as in \cite[Lemma~13]{raghavan}.
Observe that it uses neither the order type of the clubs
in the square sequence nor the cardinality of $\mathrm{FIN}$.

\begin{lemma}
There is no $X\in[\lambda]^{\omega_1}$ such that
$[X]^{\leq\omega}\subseteq\mathcal I$.
\end{lemma}

\begin{proof}
See the proof of \cite[Lemma~13]{raghavan}.
\end{proof}

From now on, assume that $\mu=\omega$.
Since $\omega<\mathfrak b$, Lemma~\ref{dilip} implies that
$\mathcal I$ is a P-ideal. Under the assumption that
$\square(\lambda,\omega)$ holds, we will show that the
second alternative of PID also fails, thereby contradicting PID.

\begin{lemma}
Suppose $X\in[\lambda]^\lambda$ satisfies
$[X]^\omega\cap\mathcal I=\emptyset$.
Then there exist $A,B\in[\lambda]^\lambda$,
$F\in\mathrm{FIN}$, and $k\in\omega$ such that
$$
\forall\alpha\in A\;\forall\beta\in B
\left(\alpha<\beta\rightarrow S_F(\alpha,\beta)\leq k\right).
$$
\end{lemma}

\begin{proof}
The proof of \cite[Lemma~14]{raghavan} applies verbatim, with
$\lambda$ in place of $\kappa^+$: it uses only the regularity
of $\lambda$, the countability of $\mathrm{FIN}$, and the
definition of $\mathcal I$, and does not use any bound on
the order types of the clubs.
\end{proof}

We are therefore  done once we prove the following lemma, which corresponds to Lemma 18 in \cite{raghavan}. However, additional work is needed here, as we are not considering the order type of the elements in square sequences.

\begin{lemma}\label{arbitrarylength}
 For every $A, B\in [\lambda]^\lambda$ and every $F\in \mathrm{FIN}$ there are  minimal walks  through $F$  starting in $B$ and finishing in $A$ of arbitrarily large length; i.e. 
 for every $k\in\omega$, there are $\alpha\in A$ and $\beta\in B$ such that $S_F(\alpha, \beta)\geq k$.  
\end{lemma}

\section*{Proof of Lemma \ref{arbitrarylength}}

Given a natural number $n\in\omega$, we equip $n^{<\omega}$, the set of all finite sequences of elements of $n$, with the \emph{shortlex order} $<_{sl}$ defined as follows: for $s,t\in n^{<\omega}$,
$$
s<_{sl}t
$$
if and only if either $|s|<|t|$, or $|s|=|t|$ and, for the least $i$ such that $s(i)\neq t(i)$, we have $s(i)<t(i)$.

\begin{definition}\label{defnested}
Let $\{C^i_\beta:(i,\beta)\in\omega\times\lambda\}$ be a
$\square(\lambda,\omega)$-sequence, and let
$$
F=\{i_l:l<|F|\}\in\mathrm{FIN}.
$$

Given $L\in\omega$, a family
$$
\{\langle\gamma_s,\delta_s,\beta_s\rangle:s\in|F|^{\le L}\}
$$
is called a \emph{system of nested intervals of length $L$} if the following conditions hold:

\begin{enumerate}
    \item $\gamma_s,\delta_s,\beta_s<\lambda$ and
    $\gamma_s<\delta_s$ for every $s\in|F|^{\le L}$;
    \item $\delta_\emptyset\le\beta_\emptyset$;
\end{enumerate}

and, letting
$$
I_s=[\gamma_s,\delta_s)
$$
for every $s\in|F|^{\le L}$,

\begin{enumerate}[resume]
    \item if $s\le_{sl}t$, then
    $
    I_s\supseteq I_t.
    $
\end{enumerate}

Moreover, for every $s\in|F|^{<L}$ and every $l<|F|$,

\begin{enumerate}[resume]
    \item
    $
    I_{s^\frown\langle l\rangle}
    \cap
    C^{i_l}_{\beta_s}
    =
    \emptyset;
    $

    \item
    $
    \beta_{s^\frown\langle l\rangle}
    =
    \min\left(
    C^{i_l}_{\beta_s}
    \setminus
    \delta_{s^\frown\langle l\rangle}
    \right).
    $
\end{enumerate}

\end{definition}

The following remark allows us to determine explicitly the initial segment
of an $F$-walk from a system of nested intervals. It will be used in the
proof of Claim~\ref{claimsequence}.

\begin{remark}\label{minim}
Let $\gamma<\delta<\beta$ be ordinals, and let $C\subseteq\beta$ with
$C\setminus\delta\neq\emptyset$. Suppose that
$$
[\gamma,\delta)\cap C=\emptyset.
$$
Then
$$
\min(C\setminus\alpha)=\min(C\setminus\delta)
$$
for every $\alpha\in[\gamma,\delta)$.
\end{remark}

\begin{proof}
Take $\alpha\in[\gamma,\delta)$. Since $\alpha\le\delta$,
$$
\min(C\setminus\alpha)
\le
\min(C\setminus\delta).
$$

On the other hand, since
$[\gamma,\delta)\cap C=\emptyset$
and
$\alpha\in[\gamma,\delta)$,
we have
$$
\min(C\setminus\alpha)\ge\delta.
$$
Therefore,
$$
\min(C\setminus\alpha)\in C\setminus\delta,
$$
and hence
$$
\min(C\setminus\alpha)
\ge
\min(C\setminus\delta).
$$

The result follows.
\end{proof}
For a given system of nested intervals $\{\langle\gamma_s, \delta_s, \beta_s\rangle : s\in |F|^{\leq L}\}$ of length $L$ with $I_s=[\gamma_s, \delta_s)$, 
let
$$
s_p=(|F|-1)^p
$$
denote the constant sequence of length $p$ with value $|F|-1$, and let
$$
I_p=I_{s_p}.
$$

\begin{claim}\label{claimsequence}
For every $p\le L$, if $\alpha\in I_p$ and
$\sigma\in T_F(\alpha,\beta_\emptyset)$,  there are
$j^\sigma_0,\ldots,j^\sigma_{p-1}<|F|$ and $i\in F$ such that
$$
\sigma=
\langle
\beta_\emptyset,
\beta_{\langle j^\sigma_0\rangle},
\ldots,
\beta_{\langle j^\sigma_0,\ldots,j^\sigma_{p-1}\rangle}
\rangle^\frown\tau,
$$
where
$$
\tau\in
T_F\!\left(
\alpha,
\min\left(
C^i_{
\beta_{\langle
j^\sigma_0,\ldots,j^\sigma_{p-1}
\rangle}
}
\setminus\alpha
\right)
\right).
$$
\end{claim}

\begin{proof}
We argue by induction on $p$.

For $p=0$, if $\alpha\in I_0$, then
$$
\alpha<\delta_\emptyset\le\beta_\emptyset.
$$
Hence, by the definition of $T_F$, every
$\sigma\in T_F(\alpha,\beta_\emptyset)$
can be written as
$$
\sigma=
\langle\beta_\emptyset\rangle^\frown\tau,
$$
where
$$
\tau\in
T_F\!\left(
\alpha,
\min(C^i_{\beta_\emptyset}\setminus\alpha)
\right)
$$
for some $i\in F$.

Suppose now that the statement holds for some $p<L$, and let
$\alpha\in I_{p+1}$ and
$\sigma\in T_F(\alpha,\beta_\emptyset)$.

Since
$s_p<_{sl}s_{p+1}$,
property (3) of Definition~\ref{defnested} gives
$$
I_{p+1}\subseteq I_p,
$$
so that $\alpha\in I_p$. By the induction hypothesis,
$$
\sigma=
\langle
\beta_\emptyset,
\beta_{\langle j^\sigma_0\rangle},
\ldots,
\beta_s
\rangle^\frown\tau,
$$
where
$$
s=
\langle
j^\sigma_0,\ldots,j^\sigma_{p-1}
\rangle,
$$
and
$$
\tau\in
T_F\!\left(
\alpha,
\min(C^i_{\beta_s}\setminus\alpha)
\right)
$$
for some $i\in F$.

Let $l<|F|$ be such that $i=i_l$.
Since
$$
\alpha\in I_{p+1}\subseteq I_{s^\frown\langle l\rangle},
$$
property (4) of Definition~\ref{defnested} yields
$$
I_{s^\frown\langle l\rangle}
\cap
C^{i_l}_{\beta_s}
=
\emptyset.
$$
Hence, by Remark~\ref{minim},
$$
\min(C^{i_l}_{\beta_s}\setminus\alpha)
=
\min\left(
C^{i_l}_{\beta_s}
\setminus
\delta_{s^\frown\langle l\rangle}
\right)
=
\beta_{s^\frown\langle l\rangle},
$$
where the last equality follows from property (5) of
Definition~\ref{defnested}. Therefore,
$$
\tau\in
T_F(\alpha,\beta_{s^\frown\langle l\rangle}).
$$

Since
$$
\alpha
<
\delta_{s_{p+1}}
\le
\delta_{s^\frown\langle l\rangle}
\le
\beta_{s^\frown\langle l\rangle},
$$
the definition of $T_F$ gives
$$
\tau=
\langle
\beta_{s^\frown\langle l\rangle}
\rangle^\frown\tau',
$$
where
$$
\tau'
\in
T_F\!\left(
\alpha,
\min(C^j_{\beta_{s^\frown\langle l\rangle}}\setminus\alpha)
\right)
$$
for some $j\in F$.

Setting
$$
j^\sigma_p=l,
$$
we obtain the desired decomposition of $\sigma$, completing the induction.
\end{proof}

\begin{corollary}\label{corolength}
If $\alpha\in I_p$, then
$$
|\sigma|\geq p+2,
$$
for every $p\leq L$ and $\sigma\in T_F(\alpha,\beta_\emptyset)$.
\end{corollary}

Thus, it remains to prove the following proposition.

\begin{proposition}\label{theproposition}
For every $\square(\lambda,\omega)$-sequence and every $F\in\mathrm{FIN}$, there exist systems of nested intervals of arbitrarily large finite length.
\end{proposition}

The following lemma is the key ingredient. It was proved in
\cite[Lemma~4.2]{torreswu2}.
\footnote{There is a typo in \cite[Lemma~4.2]{torreswu2}: it should read
``$\beta<\lambda$'' instead of ``$\beta<\mu$''.}
 
\begin{lemma}\label{lemmaunbounded}

Let $\lambda$ be a regular cardinal and $\mu$ an infinite cardinal with
$\mathrm{cof}(\mu)^+<\lambda$. Let
$\langle\mathscr C_\alpha:\alpha\in\mathrm{Lim}(\lambda)\rangle$
be a $\square(\lambda,\mu)$-sequence. For $\alpha\in\mathrm{Lim}(\lambda)$, let
$\mathscr C_\alpha=\{C^\xi_\alpha:\xi\in\mu\}$. Then, for $\theta$ sufficiently large and any
$M\prec H_\theta$ containing all relevant parameters,
with $\mu\subseteq M$ and $\sup(M\cap\lambda)<\lambda$, any
set $D\subseteq\lambda$ with
$D\cap\mathrm{cof}(>\mu)$ stationary and $D\in M$, and every $\xi\in\mu$, the set
$$
X(M,D,\xi):=_{\mathrm{def}}
\left\{
\delta\in D\cap M:
\delta\notin
\mathrm{Lim}\left(C^\xi_{\sup(M\cap\lambda)}\right)
\right\}
$$
is unbounded in $\sup(M\cap\lambda)$.

\end{lemma}

Observe that this lemma implies the existence, for a given
$M\prec H_\theta$ and every $\xi<\mu$, of an interval
$$
[\gamma(\xi),\delta(\xi))
\subseteq
\sup(M\cap\lambda)
$$
such that
$$
[\gamma(\xi),\delta(\xi))
\cap
C^\xi_{\sup(M\cap\lambda)}
=
\emptyset.
$$
For the case $\mu=\omega$, we will see that we can even obtain a system of
nested intervals in the sense of Definition~\ref{defnested}.
\begin{corollary}\label{corollaryunbounded}
Let $\lambda\geq\omega_2$ be a regular cardinal, and let
$\vec{\mathscr C}
=
\{C^i_\alpha:(i,\alpha)\in\omega\times\mathrm{Lim}(\lambda)\}$
be a $\square(\lambda,\omega)$-sequence. Then, for $\theta$ sufficiently large and any
$M\prec H_\theta$ containing all relevant parameters,
with $\sup(M\cap\lambda)<\lambda$, any club
$D\subseteq\lambda$ with $D\in M$, and every $i\in\omega$, the set
$$
\left\{
\delta\in D\cap M:
\delta\notin
\mathrm{Lim}\left(C^i_{\sup(M\cap\lambda)}\right)
\right\}
$$
is unbounded in $\sup(M\cap\lambda)$.
\end{corollary}

\begin{proof}
Since $\lambda\geq\omega_2$ is regular,
$\mathrm{cof}(>\omega)$ is stationary in $\lambda$. Therefore,
$$
D\cap\mathrm{cof}(>\omega)
$$
is stationary. 
It is a standard fact that $\omega\subseteq M$ for every elementary
submodel $M\prec H_\theta$. The conclusion follows from
Lemma~\ref{lemmaunbounded}, applied with $\mu=\omega$.
\end{proof}

\begin{proof}[Proof of Proposition~\ref{theproposition}]

Assume that the principle $\square(\lambda,\omega)$ holds, and let
$$
\vec{\mathscr C}
=
\{C^i_\beta:(i,\beta)\in\omega\times\mathrm{Lim}(\lambda)\}
$$
be a $\square(\lambda,\omega)$-sequence witnessing it.

Fix $A,B\in[\lambda]^\lambda$ and $F\in\mathrm{FIN}$.

Let $\theta>\lambda$ be a sufficiently large regular cardinal such that
all relevant parameters, as well as every club subset of $\lambda$,
belong to $H_\theta$. Let
$$
\mathfrak M_0
=
\langle H_\theta;\in,\lambda,\vec{\mathscr C},A,\ldots\rangle.
$$
Let
$$
\langle M^0_\xi:\xi<\lambda\rangle
$$
be a continuous, strictly $\subseteq$-increasing sequence of elementary
submodels of $\mathfrak M_0$, each of size less than $\lambda$, such that
$$
\bigcup_{\xi<\lambda}(M^0_\xi\cap\lambda)=\lambda.
$$
Let
$$
D_0
=
\left\{
\sup(M^0_\xi\cap\lambda):\xi<\lambda
\right\}.
$$
Then $D_0$ is a club subset of $\lambda$, and hence $D_0\in H_\theta$.

Suppose that, for some $n\in\omega$, we have constructed, for every
$l\leq n$, a continuous, strictly $\subseteq$-increasing sequence
$$
\langle M^l_\xi:\xi<\lambda\rangle
$$
of elementary submodels of $\mathfrak M_l$, each of size less than
$\lambda$, such that
$$
\bigcup_{\xi<\lambda}(M^l_\xi\cap\lambda)=\lambda,
$$
as well as the club
$$
D_l
=
\left\{
\sup(M^l_\xi\cap\lambda):\xi<\lambda
\right\}.
$$

Let
$$
\mathfrak M_{n+1}
=
\langle
H_\theta;
\in,
\lambda,
D_0,\ldots,D_n,
\vec{\mathscr C},
A,
\ldots
\rangle,
$$
and let
$$
\langle M^{n+1}_\xi:\xi<\lambda\rangle
$$
be a continuous, strictly $\subseteq$-increasing sequence of elementary
submodels of $\mathfrak M_{n+1}$, each of size less than $\lambda$, such
that
$$
\bigcup_{\xi<\lambda}(M^{n+1}_\xi\cap\lambda)=\lambda.
$$
Define
$$
D_{n+1}
=
\left\{
\sup(M^{n+1}_\xi\cap\lambda):\xi<\lambda
\right\}.
$$
Then $D_{n+1}$ is a club subset of $\lambda$, and hence
$D_{n+1}\in H_\theta$.

In particular,
$$
D_0,\ldots,D_n\in M^{n+1}_\xi
$$
for every $n\in\omega$ and every $\xi<\lambda$.

Finally, take
$$
M
\prec
\left\langle
H_\theta;
\in,
\lambda,
\vec{\mathscr C},
A,
\langle D_n:n\in\omega\rangle,
\ldots
\right\rangle
$$
with $|M|<\lambda$.

We have the following:

\begin{claim}\label{nested}
Let $F\in\mathrm{FIN}$ with
$F=\{i_l:l<|F|\}$.
Then, for every natural number $K\geq |F|$, every elementary submodel
$N\prec H_\theta$ with $\sup(N\cap\lambda)<\lambda$ such that
$D_0,D_1,\ldots,D_K\in N$ and also containing all relevant parameters,
every $\eta<\sup(N\cap\lambda)$,
and every $\beta<\lambda$ with  $\beta\geq\sup(N\cap\lambda)$, there exist a collection of
elementary submodels $\{M_l:l<|F|\}$ of $H_\theta$ and ordinals
$\{\gamma_l:l<|F|\}$ such that

\begin{enumerate}
    \item $\gamma_0\geq\eta$;

    \item
    $\sup(M_l\cap\lambda)\in D_{K-l}$;

    \item
    $\gamma_l\in M_l\cap \lambda$ and $ D_0,\ldots,D_{K-l-1}\in M_l$;

    \item
    $[\gamma_l,\sup(M_l\cap\lambda))
    \cap C^{i_l}_\beta=\emptyset$;

    \item
    $\sup(M_0\cap\lambda)<\sup(N\cap\lambda)$;

    \item
    $\sup(M_{l+1}\cap\lambda)<\sup(M_l\cap\lambda)$
    for every $l<|F|-1$;

    \item
    $[\gamma_l,\sup(M_l\cap\lambda))
    \supseteq
    [\gamma_{l+1},\sup(M_{l+1}\cap\lambda))$
    for every $l<|F|-1$.
\end{enumerate}
\end{claim}

\begin{proof}
We proceed by induction on $l<|F|$.

First let $l=0$. Put $\delta_N=\sup(N\cap\lambda)$.

If $\beta=\delta_N$, apply
Corollary~\ref{corollaryunbounded} directly to
$N$, $D_K$, and $i_0$ to find $\xi<\lambda$ such that
$$
\delta_0:=\sup(M^K_\xi\cap\lambda)
\in
D_K\cap N
\setminus
\left(
\mathrm{Lim}(C^{i_0}_\beta)
\cup(\eta+1)
\right).
$$
Let $M_0=M^K_\xi$.
Since $\delta_0\notin\mathrm{Lim}(C^{i_0}_\beta)$,
there is $\gamma<\delta_0$ such that
$$
[\gamma,\delta_0)\cap C^{i_0}_\beta=\emptyset.
$$
Choose $\gamma_0\in M_0\cap\lambda$ with
$\gamma_0\geq\max\{\gamma,\eta\}$.
Then
$$
[\gamma_0,\delta_0)\cap C^{i_0}_\beta=\emptyset.
$$

Now suppose that $\beta>\delta_N$.
We distinguish two cases.

\medskip

\noindent
\underline{Case 1.}
$\delta_N\notin\mathrm{Lim}(C^{i_0}_\beta)$.

In particular, there is $\gamma<\delta_N$ such that
$$
[\gamma,\delta_N)\cap C^{i_0}_\beta=\emptyset.
$$
Since $D_K\in N$, $D_K\cap N$ is unbounded in $\delta_N$.
Choose $\xi<\lambda$ such that
$$
\sup(M^K_\xi\cap\lambda)
\in
D_K\cap N\setminus(\max\{\gamma,\eta\}+1).
$$
Then choose
$\gamma'\geq\max\{\gamma,\eta\}$
with $\gamma'\in M^K_\xi\cap\lambda$.
Set $\gamma_0=\gamma'$ and $M_0=M^K_\xi$.
Thus
$$
[\gamma_0,\sup(M_0\cap\lambda))
\cap C^{i_0}_\beta=\emptyset.
$$

\medskip

\noindent
\underline{Case 2.}
$\delta_N\in\mathrm{Lim}(C^{i_0}_\beta)$.

By coherence, there is $j\in\omega$ such that
$$
C^{i_0}_\beta\cap\delta_N=C^j_{\delta_N}.
$$

Apply Corollary~\ref{corollaryunbounded} to
$N$, $D_K$, and $j$ to find $\xi<\lambda$ such that
$$
\sup(M^K_\xi\cap\lambda)
\in
D_K\cap N
\setminus
\left(
\mathrm{Lim}(C^j_{\delta_N})
\cup(\eta+1)
\right).
$$

Let $M_0=M^K_\xi$, and choose
$\gamma_0\in M_0\cap\lambda$ with $\gamma_0\geq\eta$
such that
$$
[\gamma_0,\sup(M_0\cap\lambda))
\cap C^j_{\delta_N}=\emptyset.
$$

Since
$$
C^j_{\delta_N}=C^{i_0}_\beta\cap\delta_N,
$$
it follows that
$$
[\gamma_0,\sup(M_0\cap\lambda))
\cap C^{i_0}_\beta=\emptyset.
$$

In all three cases,
$D_0,\ldots,D_{K-1}\in M_0$.

Now suppose that $l+1<|F|$ and that we have already
constructed a submodel $M_l$ and an ordinal $\gamma_l$
such that

\begin{enumerate}
    \item
    $\sup(M_l\cap\lambda)\in D_{K-l}$;

    \item
    $\gamma_l\in M_l\cap \lambda$ and $D_0,\ldots,D_{K-l-1}\in M_l$;

    \item
    $[\gamma_l,\sup(M_l\cap\lambda))
    \cap C^{i_l}_\beta=\emptyset$.
\end{enumerate}

We next construct a submodel $M_{l+1}$ and an ordinal
$\gamma_{l+1}$ satisfying

\begin{enumerate}
    \item
    $\sup(M_{l+1}\cap\lambda)\in D_{K-l-1}$;

    \item
    $\gamma_{l+1}\in M_{l+1}\cap \lambda$ and $D_0,\ldots,D_{K-l-2}\in M_{l+1}$;

    \item
    $[\gamma_{l+1},\sup(M_{l+1}\cap\lambda))
    \cap C^{i_{l+1}}_\beta=\emptyset$;

    \item
    $\sup(M_{l+1}\cap\lambda)
    <\sup(M_l\cap\lambda)$;

    \item
    $[\gamma_{l+1},\sup(M_{l+1}\cap\lambda))
    \subseteq
    [\gamma_l,\sup(M_l\cap\lambda))$.
\end{enumerate}
Put $\delta_l=\sup(M_l\cap\lambda)$.
By construction,
$\delta_l\leq\delta_0<\delta_N\leq\beta$.
Again we distinguish two cases.
\medskip

\noindent
\underline{Case 1.}
$\delta_l\notin\mathrm{Lim}(C^{i_{l+1}}_\beta)$.

In particular, there is $\gamma<\delta_l$ such that
$$
[\gamma,\delta_l)\cap C^{i_{l+1}}_\beta=\emptyset.
$$
Let $\gamma^*=\max\{\gamma,\gamma_l\}$.

Since $D_{K-l-1}\in M_l$,
$D_{K-l-1}\cap M_l$ is unbounded in $\delta_l$.
Choose $\xi<\lambda$ such that
$$
\sup(M^{K-l-1}_\xi\cap\lambda)
\in
D_{K-l-1}\cap M_l\setminus(\gamma^*+1),
$$
and choose $\gamma'\geq\gamma^*$ with
$\gamma'\in M^{K-l-1}_\xi\cap\lambda$.

Set
$\gamma_{l+1}=\gamma'$ and
$M_{l+1}=M^{K-l-1}_\xi$.
Thus
$$
[\gamma_{l+1},\sup(M_{l+1}\cap\lambda))
\cap C^{i_{l+1}}_\beta=\emptyset.
$$

\medskip

\noindent
\underline{Case 2.}
$\delta_l\in\mathrm{Lim}(C^{i_{l+1}}_\beta)$.

By coherence, there is $j\in\omega$ such that
$$
C^{i_{l+1}}_\beta\cap\delta_l=C^j_{\delta_l}.
$$

Apply Corollary~\ref{corollaryunbounded} to
$M_l$, $D_{K-l-1}$, and $j$ to find
$\xi<\lambda$ such that
$$
\sup(M^{K-l-1}_\xi\cap\lambda)
\in
D_{K-l-1}\cap M_l\setminus(\gamma_l+1),
$$
and
$$
\sup(M^{K-l-1}_\xi\cap\lambda)
\notin
\mathrm{Lim}(C^j_{\delta_l}).
$$

Let $M_{l+1}=M^{K-l-1}_\xi$,
and choose
$\gamma_{l+1}\in M_{l+1}\cap\lambda$
with $\gamma_{l+1}\geq\gamma_l$ such that
$$
[\gamma_{l+1},\sup(M_{l+1}\cap\lambda))
\cap C^j_{\delta_l}=\emptyset.
$$

Since
$$
C^j_{\delta_l}=C^{i_{l+1}}_\beta\cap\delta_l,
$$
it follows that
$$
[\gamma_{l+1},\sup(M_{l+1}\cap\lambda))
\cap C^{i_{l+1}}_\beta=\emptyset.
$$

In both cases,
$D_0,\ldots,D_{K-l-2}\in M_{l+1}$.

This completes the induction and hence the proof of the claim.
\end{proof}

\newcommand{\opl}
    {
            {
                {\{0\}^p
                }^\frown\langle l \rangle
            }
    }

\newcommand{\oplo}
    {
            {
                {\{0\}^p
                }^\frown\langle l+1 \rangle
            }
    }

Fix a natural number $k$ and let $n=|F|$.
Choose $K\in\omega$ sufficiently large,  at least the cardinality of $n^{\leq k}$ (i.e. $K\geq\sum\limits_{i=0}^k n^i$), 
 so that all applications of
Claim~\ref{nested} can be carried out using  the club sets
$D_0,\ldots,D_K. $

We construct by induction a system of nested intervals
$$
\left\{
\langle\gamma_s,\delta_s,\beta_s\rangle:
s\in n^{\leq k}
\right\},
$$
together with elementary submodels
$$
\{M_s:s\in n^{\leq k}\},
$$
such that
$$
\delta_s=\sup(M_s\cap\lambda)
$$
for every $s\in n^{\leq k}$.

Let
$
M_\emptyset=M,
$
$
\delta_\emptyset=\sup(M\cap\lambda)
$
and
$
\gamma_\emptyset= \emptyset.
$
Choose $\beta \in B$ and 
a limit ordinal
$
\beta_\emptyset>\delta_\emptyset
$
such that $\beta=  \beta_\emptyset+m$ for some $m\in \omega$.  This choice is possible since $B$ is unbounded in $\lambda$. 
Observe also that for every $\alpha<\beta_\emptyset$, we have   $S_F(\alpha, \beta_\emptyset)\leq S_F(\alpha, \beta) $.  Properties (1) and (2) of Definition \ref{defnested}
hold at the root. 

Apply Claim~\ref{nested} to $M_\emptyset$,  $\beta_\emptyset$, $\gamma_\emptyset$ and $K$ to obtain
a collection of submodels
$\{M_{\langle l\rangle}:l<n\}$
and ordinals
$\{\gamma_{\langle l\rangle}:l<n\}$
such that

\begin{enumerate}
    \item  $\sup (M_{\langle 0 \rangle }\cap \lambda)<\sup (M_{\emptyset}\cap \lambda) $; 
    
    \item $\sup (M_{\langle l \rangle}\cap \lambda)\in D_{K-l}$; 
    
    \item
    $\gamma_{\langle l\rangle}\in M_{\langle l\rangle}\cap \lambda$ and $
    D_0,\ldots,D_{K-l-1}
    \in M_{\langle l\rangle}$;

    \item
    $[\gamma_{\langle l\rangle},
    \sup(M_{\langle l\rangle}\cap\lambda))
    \cap
    C^{i_l}_{\beta_\emptyset}
    =
    \emptyset$;

    \item 
    $\sup(M_{\langle l\rangle}\cap\lambda)
    >
    \sup(M_{\langle l+1\rangle}\cap\lambda)$
    for every $l<n-1$;
    
    \item
    $[\gamma_{\langle l\rangle},
    \sup(M_{\langle l\rangle}\cap\lambda))
    \supseteq
    [\gamma_{\langle l+1\rangle},
    \sup(M_{\langle l+1\rangle}\cap\lambda))$
    for every $l<n-1$.
\end{enumerate}

For every $l<n$, let
$$
\delta_{\langle l\rangle}
=
\sup(M_{\langle l\rangle}\cap\lambda),
$$
and define
$$
\beta_{\langle l\rangle}
=
\min\left(
C^{i_l}_{\beta_\emptyset}
\setminus
\delta_{\langle l\rangle}
\right).
$$

Observe in particular  that $\beta_{\langle l\rangle }\geq \delta_{\langle  l \rangle}$ for all $l<n$ and  $\delta_{\langle  l \rangle}>\delta_{\langle n-1\rangle}$ for $l<n-1$. All five conditions of Definition~\ref{defnested}
hold on $n^{\leq 1}$. 

Suppose that, for some $p\geq 1$, we have already constructed, for every $s\in n^{\leq p}$, triples
$\langle\gamma_s,\delta_s,\beta_s\rangle$
and elementary submodels
$M_s$
such that
$$
\delta_s=\sup(M_s\cap\lambda),
$$
where  $\beta_s= \min \left( C^{i_l}_{\beta_t}\setminus \delta_s \right) $ whenever $s= t^\frown \langle  l \rangle $, 
and satisfying the following properties:

\begin{enumerate}

\item
$\delta_s>\delta_{s'}$
whenever $s<_{sl}s_p$;

\item
$\delta_s\in D_{K-\varphi(s)}$
for every nonempty $s\in n^{\leq p}$;

\item
$\gamma_s\in  M_s\cap \lambda$  and $
D_0,\ldots, D_{K-\varphi(s')}
\in M_s$ for every nonempty $s\in n^{\leq p}$;

\item
$I_s\cap C^{i_l}_{\beta_t}=\emptyset$
whenever
$s=t^\frown\langle l\rangle$;

\item
$I_s\supseteq I_{s'}$
whenever $s<_{sl}s_p$.

\end{enumerate}

Here $s'$ denotes the $<_{sl}$-successor of $s$, and

$$
\varphi:
\left(n^{<\omega}\setminus\{\emptyset\},<_{sl}\right)
\longrightarrow
(\omega,<)
$$
is the canonical order-preserving bijection.

Observe that $\beta_s\geq \delta_s$ for every $s\in n^{\leq p}$ and in particular, if $s=t^\frown \langle l \rangle $, we have   $\beta_{t^\frown\langle l\rangle }>\delta_{t^\frown \langle n-1\rangle }$ for $l<n-1$ and properties (1)-(5) of a system of nested intervals hold for $s\in n^{\leq p}$. 

We now construct the triples
$\langle\gamma_t,\delta_t,\beta_t\rangle$ and elementary submodels $M_t$ with  $\delta_t=\sup (M_t\cap \lambda)$ for $t\in n^{p+1}$.
We first construct the intervals
$I_t=[\gamma_t,\delta_t)$ and submodels $M_t$ for $t\in n^{p+1}$ using the
$\beta_s$'s for $s\in n^p$. At the end, we define the corresponding
$\beta_t$'s for $t\in n^{p+1}$.

We proceed by induction on $n^p$ with respect to the order $<_{sl}$ and
for each $s\in n^p$, we construct the intervals
$I_{s^\frown\langle l\rangle}$ 
and submodels $M_{s\frown\langle l \rangle}$
for $l<n$ using Claim \ref{nested}.

Consider $\beta_{\{0\}^p}$ and observe that $D_0, \ldots, D_{K-\varphi\left( \{0\}^{p+1}\right)}\in M_{s_p}$.  
Apply Claim~\ref{nested} to $\beta_{\{0\}^p}$,  $M_{s_p}$,  $\gamma_{s_p}$  and 
$K-\varphi\left(\{0\}^{p+1} \right) $ 
to build 
$\gamma_{{\{0\}^p}^\frown \langle l \rangle}$ and 
$M_{{\{0\}^p}^\frown \langle l \rangle }$ such that: 

\begin{enumerate}
    \item $\gamma_{\{0\}^{p+1}}\geq \gamma_{s_p}$;
    \item $\delta_{\{0\}^{p+1}}<\delta_{s_p}$; 
    \item $\delta_\oplo<\delta_\opl$ for all $l<n-1$;
    \item $\delta_\opl\in D_{K-\varphi \left(\opl\right )}$ for all $l<n$;
    \item $\gamma_\opl\in M_\opl\cap \lambda$ and $ D_0, \ldots, D_{K-\varphi\left(\opl\right)-1} \in M_\opl$ for all $l<n$;
    \item $I_\opl\cap C^{i_l}_{\beta_{\{0\}^p}}= \emptyset$ for all $l<n$; 
    \item $I_\oplo\subseteq I_\opl$ for all $l<n-1$. 
\end{enumerate}

Now suppose that the families
$$
\left\{ \left(
\gamma_{s^\frown\langle l\rangle}, \delta_{s^\frown\langle l\rangle}\right) :l<n
\right\}
$$
have already been constructed together with the submodels $\{M_{s^\frown\langle l  \rangle }: l<n\}$, where $\delta_{s^\frown \langle l\rangle }= \sup (M_{s^\frown \langle l\rangle }\cap \lambda)$  for $s\in n^p$ and $s<_{sl}s_p$ such that:  

\newcommand{\spl}
    {   
        {
        s'^\frown\langle l \rangle  
        }
    }

\newcommand{\splo}
    {
        {
    s'^\frown\langle l+1 \rangle
        }
    }

\newcommand{\esl}
    {   
        {
        s^\frown\langle l \rangle  
        }
    }

\newcommand{\eslo}
    {
        {
    s^\frown\langle l+1 \rangle
        }
    }

\begin{enumerate}
    \item $\gamma_{s^\frown\langle 0\rangle}\geq \gamma_t$ 
    and $\delta_{s^\frown\langle 0\rangle}<\delta_t$
    for every $t<_{sl} s^\frown \langle 0 \rangle$; 
    \item $\delta_\eslo<\delta_\esl$ for all $l<n-1$;
    \item $\delta_\esl\in D_{K-\varphi \left(\esl\right )}$ for all $l<n$;
    \item $\gamma_\esl\in M_\esl \cap \lambda$ and $ D_0, \ldots,D_{K-\varphi\left(\esl\right)-1} \in M_\esl$ for all $l<n$;
    \item $I_\esl\cap C^{i_l}_{\beta_{s}}= \emptyset$ for all $l<n$; 
    \item $I_\eslo\subseteq I_\esl$ for all $l<n-1$. 
\end{enumerate}

Let $s'$ be the $<_{sl}$-successor of $s$.   

Apply Claim \ref{nested} to $M_{s^\frown \langle n-1\rangle}$, $\beta_{s'}$, $\gamma_{s^\frown\langle n-1 \rangle }$ and  $K-\varphi ({s'}^\frown \langle  0\rangle )$ to build the respective 
 families
$$
\left\{ \left(
\gamma_{s'^\frown\langle l\rangle}, \delta_{s'^\frown\langle l\rangle}\right) :l<n
\right\}
$$
together with the submodels $\{M_{s'^\frown\langle l  \rangle }: l<n\}$, where 
$$\delta_{s'^\frown \langle l\rangle }= \sup (M_{s'^\frown \langle l\rangle }\cap \lambda),$$ 
and such that: 

\begin{enumerate}
    \item $\gamma_{s'^\frown\langle 0\rangle}\geq \gamma_{s^\frown\langle n-1\rangle }$;
    \item $\delta_{s'^\frown\langle 0\rangle}<\delta_{s^\frown\langle n-1\rangle }$; 
    \item $\delta_\splo<\delta_\spl$ for all $l<n-1$;
    \item $\delta_\spl\in D_{K-\varphi \left(\spl\right )}$ for all $l<n$;
    \item $\gamma_\spl\in M_\spl\cap \lambda$ and $ D_0, \ldots, D_{K-\varphi\left(\spl\right)-1} \in M_\spl$ for all $l<n$;
    \item $I_\spl\cap C^{i_l}_{\beta_{s'}}= \emptyset$ for all $l<n$; 
    \item $I_\splo\subseteq I_\spl$ for all $l<n-1$. 
\end{enumerate}

This completes the  induction over $n^p$ with respect to the order  $<_{sl}$ and the construction of the corresponding families 
$$
\left\{ \left(
\gamma_{s^\frown\langle l\rangle}, \delta_{s^\frown\langle l\rangle}\right) :l<n
\right\}
$$
together with the submodels $\{M_{s^\frown\langle l  \rangle }: l<n\}$, where 
$$\delta_{s^\frown \langle l\rangle }= \sup (M_{s^\frown \langle l\rangle }\cap \lambda)$$  
for  all $s\in n^p$.   For $s\in n^p$, let  
$$\beta_{s^\frown \langle l \rangle}= \min \left(C^{i_l}_{\beta_s}\setminus \delta_{s^\frown \langle l \rangle }\right)  $$
for $l<n$. 
Observe in particular that $\beta_{s^\frown \langle l \rangle}\geq \delta_{s^\frown \langle l \rangle}>\delta_{s^\frown \langle n-1 \rangle} $ for $l<n-1$ and properties (1)-(5)  of a system of nested intervals hold for $s\in n^{\leq p+1}$. 

With this,  we have already built the families 
$$
\{\left(\gamma_s, \delta_s, \beta_s, M_s\right): s\in n^{p+1}\}, 
$$
and this finishes the induction over the length of the sequences. 

In particular, 
$\{\langle \gamma_s, \delta_s, \beta_s\rangle: s\in n^{\leq k} \}$ is a system of length $k$. 

This proves Proposition~\ref{theproposition}.

\end{proof}

We now deduce Lemma~\ref{arbitrarylength} from the construction.
Observe that $A\in M_s$ for every $s\in n^{\leq k}$. 
In particular, we can take  $\alpha\in I_k\cap A$. Apply Corollary \ref{corolength} to $\alpha$ and $\beta_\emptyset$ and so $S_F(\alpha, \beta_\emptyset)\geq k$. Since $S_F(\alpha, \beta)\geq S_F(\alpha, \beta_\emptyset)\geq k$, this completes the proof of Lemma~\ref{arbitrarylength} and the proof of Theorem \ref{pidnonsquare}.

\bibliographystyle{plain}
 \bibliography{torreswupidbib}

\end{document}